\documentclass{amsart}
\usepackage{amssymb}
\usepackage{amsmath}

\usepackage{diagrams}
\diagramstyle[centredisplay,dpi=600,UglyObsolete,heads=LaTeX]

\newcommand\Zset{\mathbb {Z}}

\newcommand\Qmax{Q^r_{\mathrm{max}}}

\newcommand\Qsimmax{Q^{\sigma}_{\mathrm{max}}}
\newcommand\Qtot{Q^r_{\mathrm{tot}}}

\newcommand\Qsimtot{Q^{\sigma}_{\mathrm{tot}}}

\newcommand\Qcl{Q^r_{\mathrm{cl}}}
\newcommand\Qlrcl{Q_{\mathrm{cl}}}
\newcommand\f{{\mathcal F}}
\newcommand\te{{\mathcal T}}
\newcommand\ef{{\mathfrak F}}

\newcommand\dirlim{\mathop{\varinjlim}\limits}
\newcommand\homo{\mathrm{Hom}}
\newcommand\ann{{\mathrm{ann}}}

\newtheorem{theorem}{Theorem}
\newtheorem{lemma}{Lemma}
\newtheorem{corollary}{Corollary}
\newtheorem{prop}{Proposition}

\newtheorem{remark}{Remark}
\newtheorem{definition}{Definition}

\begin{document}

\title{Extending Higher Derivations to Rings and Modules of Quotients}

\author{Lia Va\v s}
\address{Department of Mathematics, Physics and Statistics\\ University of the Sciences
in Philadelphia\\ Philadelphia, PA 19104, USA}
\email{l.vas@usip.edu}



\maketitle

\begin{abstract} A torsion theory is called differential (higher differential) if a derivation (higher derivation) can be extended from any module to the module of quotients corresponding to the torsion theory. We study conditions equivalent to higher differentiability of a torsion theory. It is known that the Lambek, Goldie and any perfect torsion theories are differential. We show that these classes of torsion theories are  higher differential as well. Then, we study conditions under which a higher derivation extended to a right module of quotients extends also to a right module of quotients with respect to a larger torsion theory. Lastly, we define and study the symmetric version of higher differential torsion theories. We prove that the symmetric versions of the results on higher differential (one-sided) torsion theories hold for higher derivations on symmetric modules of quotients. In particular, we prove that the symmetric Lambek, Goldie and any perfect torsion theories are higher differential.
\end{abstract}

\section{Introduction}
\label{section_introduction}

A a {\em derivation} on a ring $R$ is a mapping $\delta: R \rightarrow R$ such that $\delta(r+s)=\delta(r)+\delta(s)$ and $\delta(rs)=\delta(r)s+r\delta(s)$ for all $r,s\in R.$ A mapping $d: M\rightarrow M$ on a right $R$-module $M$ is a {\em $\delta$-derivation} if $d(x+y)=d(x)+d(y)$ and $d(xr)=d(x)r+x\delta(r)$ for all $x\in M$ and $r\in R.$  In \cite{Golan_paper}, \cite{Bland_paper}, \cite{Lia_Diff} and \cite{Lia_Extending}, the authors study how derivations agree with an arbitrary hereditary torsion theory for that ring and the conditions under which one can extend a derivation from a ring (module) to the ring of quotients (module of quotients). In particular, in \cite{Lia_Diff} and \cite{Lia_Extending} it is shown that such extension is possible for some important classes of rings and modules of quotients.

Torsion theories provide a basis for a uniform treatment of different rings of quotients and a framework suitable for a comprehensive study of rings and modules of quotients. We shall use torsion theories extensively throughout this paper. We use the usual definition of torsion theory and hereditary torsion theory (e.g. \cite{Stenstrom}, \cite{Bland_book}, \cite{Bland_paper}, \cite{Lia_Diff}). If $\tau=(\te, \f)$ is a torsion theory for $R,$ we denote the torsion submodule of a right $R$-module $M$ by $\te(M).$ If $\te(R)=0,$ $\tau$ is said to be faithful. If $\tau$ is hereditary, we denote its Gabriel filter (the collection of all right ideals $I$ such $R/I$ is a torsion module) by $\ef.$

If $\tau$ is a hereditary torsion theory with Gabriel filter $\ef$ and $M$ is a right $R$-module, the module of quotients $M_{\ef}$ of $M$ is defined as the largest submodule $N$ of the injective envelope $E(M/\te(M))$ of $M/\te(M)$ such that $N/(M/\te(M))$ is torsion module (i.e. the closure
of $M/\te(M)$ in $E(M/\te(M))$). Equivalently, the module of quotients $M_{\ef}$ can be defined by
\[M_{\ef}=\dirlim_{I\in\ef}\homo(I, \frac{M}{\te(M)})\]
(see chapter IX in \cite{Stenstrom}). Note that from this description it directly follows that $M_{\ef}=(M/\te(M))_{\ef}.$

The $R$-module $R_{\ef}$ has a ring structure and $M_{\ef}$ has a structure of a right
$R_{\ef}$-module (see exposition on pages 195--197 in \cite{Stenstrom}). The ring $R_{\ef}$ is called the
right ring of quotients with respect to the torsion theory $\tau.$

Consider the map $q_M:M\rightarrow M_{\ef}$ obtained by
composing the projection $M\rightarrow M/\te(M)$ with the injection $M/\te(M)\rightarrow M_{\ef}.$ Note that $q_M: m\mapsto L_m$ where $L_m$ is the left multiplication by $m.$ The kernel and cokernel of $q_M$ are torsion modules and $M_{\ef}$ is torsion-free (Lemmas 1.2 and 1.5, page 196, in \cite{Stenstrom}). The map $M \mapsto q_M$ defines a left exact functor $q$ from the category of right $R$-modules to the category of right $R_{\ef}$-modules (see \cite{Stenstrom} pages 197--199). The functor $q$ maps torsion modules to 0 (see Lemma 1.3, page 196 in  \cite{Stenstrom}).

A Gabriel filter $\ef$ is a {\em differential filter} if for every $I\in \ef$ there is $J\in \ef$ such that $\delta(J)\subseteq I$ for all derivations $\delta.$ The hereditary torsion theory determined by $\ef$ is said to be {\em differential} in this case. By Lemma 1.5 from \cite{Bland_paper}, a torsion theory is differential if and only if
\[d(\te(M))\subseteq \te(M)\]
for every right $R$-module $M,$ every derivation $\delta$ and every $\delta$-derivation $d$ on $M.$

In Theorem on page 277 and Corollary 1 on page 279 of \cite{Golan_paper}, Golan has shown the following
\begin{prop}[Golan] Let $\delta$ be a derivation on $R$, $M$ a right $R$-module, $d$ a $\delta$-derivation on $M$ and $\tau=(\te, \f)$ a hereditary torsion theory.
\begin{enumerate}
\item If $M$ is torsion-free, then $d$ extends to a derivation on the module of quotients $M_{\ef}$ such that  $d q_M=q_M d.$

\item If $d(\te(M))\subseteq \te(M),$ then $d$ extends to a derivation on the module of quotients $M_{\ef}$ such that $d q_M=q_M d.$
\end{enumerate}
\label{Golan_Proposition}
\end{prop}

A direct corollary of the first part of Proposition \ref{Golan_Proposition} is that a ring derivation extends to a right ring of quotients for every hereditary and faithful torsion theory. By the second part of Proposition \ref{Golan_Proposition}, we can extend a derivation on a module to a derivation on its module of quotients for every differential torsion theory. Bland proved that such extension is unique and that the converse is also true. Namely, Propositions 2.1 and 2.3 of his paper \cite{Bland_paper} state the following.

\begin{theorem}[Bland] Let $\ef$ be a Gabriel filter.
\begin{enumerate}
\item If a $\delta$-derivation on a module $M$ extends to a $\delta$-derivation on the module of quotients $M_{\ef},$ then such extension is unique.

\item The filter $\ef$ is differential if and only if for every ring derivation $\delta,$
every $\delta$-derivation on any module $M$ extends uniquely to a $\delta$-derivation on the module of quotients $M_{\ef}.$
\end{enumerate}
\label{Bland_Theorem}
\end{theorem}

In \cite{Rim2} and \cite{Bland_paperHD}, it has been studied how {\em higher} derivations agree with torsion theories. This has lead to conditions under which one can extend a higher derivation from a ring (module) to the ring of quotients (module of quotients).

Recall that a {\em higher derivation (HD) of order $n$} on $R$ is an indexed family $\{\delta_i\}_{i=0}^n$ of additive maps $\delta_i$ such that $\delta_0$ is the identity mapping on $R$ and
\[\delta_i(rs)=\sum_{j=0}^i \delta_j(r)\delta_{i-j}(s)\]
for all $i\leq n.$ An indexed family $\Delta=\{\delta_n\}_{n\in\omega}$ is a {\em higher derivation (HD)} on $R$ if $\{\delta_i\}_{i=0}^n$ is a higher derivation of order $n$ for all $n\in\omega.$

Let $\Delta$ be a higher derivation on $R.$ If $\{d_i\}_{i=0}^n$ is an indexed family of additive maps on a right $R$-module $M$ such that $d_0$ is the identity mapping on $M$ and $d_i(mr)=\sum_{j=0}^i d_j(m)\delta_{i-j}(r)$ for all $i\leq n,$ we say that $\{d_i\}$ is {\em higher $\Delta$-derivation of order $n$} on $M.$ An indexed family $D=\{d_n\}_{n\in\omega}$ is a higher $\Delta$-derivation
($\Delta$-HD for short) on $M$ if $\{d_i\}_{i=0}^n$ is a higher $\Delta$-derivation of order $n$ for all $n\in\omega.$
Occasionally we shall be adding the superscript $^M$ to $\{d_n\}$ and use $\{d^M_n\}$ to emphasize that $d_n$ is defined on a module $M.$

For example, if $\delta$ is a derivation on $R,$ then $\Delta=\{\frac{\delta^i}{i!}\}$ is a higher derivation on $R$ and if $d$ is any $\delta$-derivation on a module $M$, then $\{\frac{d^i}{i!}\}$ is a higher $\Delta$-derivation on $M$.

A Gabriel filter $\ef$ is a {\em higher differential filter of order $n$} if for every HD $\{\delta_i\}_{i=0}^n$ of order $n,$ and every $I\in \ef,$ there is $J\in \ef$ such that $\delta_i(J)\subseteq I$ for all $i\leq n.$ The hereditary torsion theory determined by $\ef$ is said to be {\em higher differential of order $n$} in this case. Since a finite intersection of right ideals from $\ef$ is again in $\ef,$ $\ef$ is a HD filter of order $n$ iff for every HD $\{\delta_i\}_{i=0}^n$ of order $n$, every $I\in \ef,$ and every $i\leq n,$ there is $J\in \ef$ such that $\delta_i(J)\subseteq I.$

A Gabriel filter $\ef$ is a {\em higher differential filter} (HD filter for short) if $\ef$ is higher differential of order $n$ for every $n.$ The hereditary torsion theory determined by $\ef$ is said to be {\em higher differential} (HD for short) in this case. Since a finite intersection of right ideals from $\ef$ is again in $\ef,$ $\ef$ is a HD filter iff for every HD $\{\delta_n\},$ every $I\in \ef$ and every $n,$ there is $J\in \ef$ such that $\delta_n(J)\subseteq I.$
Example 3.2 in \cite{Bland_paperHD} demonstrates that every hereditary torsion theory of a commutative ring is higher differential.

Let $\Delta=\{\delta_n\}$ be a HD on $R$, $M$ a right $R$-module, $D=\{d_n\}$ a $\Delta$-HD on $M$ and $\tau$ a hereditary torsion theory. If $D$ is such that every $d_n$ extends to the module of quotients $M_{\ef}$ such that $d_n q_M=q_M d_n$ for all $n,$ then we shall say that $D$ {\em extends} to a $\Delta$-HD on $M_{\ef}.$

In section \ref{section_HD}, we review Rim's and Bland's results (from \cite{Rim2} and \cite{Bland_paperHD} respectively) on higher derivations that generalize Golan's and Bland's results on differentiability of torsion theories and the conditions under which one can extend a derivation from a module to the module of quotients. In particular, we study conditions equivalent to  higher differentiability of a torsion theory.

In section \ref{section_Lambek_Goldie_perfect}, we prove that some important torsion theories are HD (Theorem \ref{3TT_are_HD}). First, we consider the Lambek torsion theory -- the torsion theory cogenerated by the injective envelope $E(R)$ (see sections 13B and 13C in \cite{Lam} and Example 1, page 200, \cite{Stenstrom}). The Lambek torsion theory is the largest faithful hereditary torsion theory. The Gabriel filter of this torsion theory is the set of all dense right ideals (see definition 8.2. in \cite{Lam} and Proposition VI 5.5, p. 147 in \cite{Stenstrom}). The ring of quotients with respect to the Lambek torsion theory is the maximal right ring of quotients $\Qmax(R).$ By Proposition 9 in \cite{Lia_Diff}, the Lambek torsion theory is differential. In the first part of Theorem \ref{3TT_are_HD}, we show that the Lambek torsion theory is higher differential.

Secondly, we consider the Goldie torsion theory. This hereditary torsion theory is defined by the condition that its torsion-free class is the class of nonsingular modules. It is larger than any hereditary faithful torsion theory (see Example 3, p. 26 in \cite{Bland_book}). So, the Lambek torsion theory is smaller than the Goldie's. The Goldie torsion theory is differential (Proposition 14 in \cite{Lia_Diff}). In the second part of Theorem \ref{3TT_are_HD} we show that the Goldie torsion theory is HD.

Lastly, we consider a class of torsion theories that represents a generalization of the classical torsion theory. Recall that the classical torsion theory is defined for any right Ore ring by the condition that the set of torsion elements of a module $M$ is $\{m\in M\,|\,mr =0$ for some regular element $r$ of $R\}.$ It is hereditary and faithful. Its right ring of quotients is called the classical right ring of quotients and is denoted by $\Qcl(R).$ The classical torsion theory has the property that $M_{\ef}\cong M\otimes_R \Qcl(R).$ The class of perfect Gabriel filters generalizes this property. If a Gabriel filter $\ef$ has the property \[M_{\ef}\cong M\otimes_R R_{\ef}\] for every right $R$-module $M,$ then $\ef$ is called perfect and the right ring of quotients $R_{\ef}$ is called the perfect right ring of quotients. These filters are interesting as all modules of quotients are determined solely by the right ring of quotients.  Perfect right rings of quotients and perfect filters have been studied and characterized in more details (see Theorem 2.1, page 227 in \cite{Stenstrom} and Proposition 3.4, page 231 in \cite{Stenstrom}). By Proposition 4 in \cite{Lia_Diff}, every perfect filter is differential. In the third part of Theorem \ref{3TT_are_HD}, we show that every perfect torsion theory is HD.

Throughout sections \ref{section_HD} and \ref{section_Lambek_Goldie_perfect}, we shall use $(r:I)$ to denote the set $\{s\in R\,|\,rs\in I\}$ where $I$ is a right ideal and $r$ an element of $R.$

In section \ref{section_agreement}, we study conditions under which extensions of a HD on two different modules of quotients agree. Let $\ef_1$ and $\ef_2$ be two filters such that $\ef_1\subseteq \ef_2.$ If $M$ is a right $R$-module and $q_i$ the natural maps $M\rightarrow M_{\ef_i}$ for $i=1,2,$ let $q_{12}$ denote the map $M_{\ef_1}\rightarrow M_{\ef_2}$ induced by inclusion  $\ef_1\subseteq \ef_2.$ In this case $q_{12}q_1=q_2.$ Let $d$ be a derivation on $M$ that extends to $M_{\ef_1}.$ Conditions under which $d$ can be extended to $M_{\ef_2}$ so that the following diagram commutes were studied in \cite{Lia_Extending}.
\begin{diagram}
  &             & M_{\ef_1}    &      & \rTo^{d} &             & M_{\ef_1}\\
  & \ruTo^{q_1} & \vLine       &      &            & \ruTo^{q_1} &    \\
M &             & \rTo^{d}     &      & M          &             & \dTo_{q_{12}}\\
  & \rdTo^{q_2} & \dTo_{q_{12}}&      &            & \rdTo^{q_2} &     \\
  &             & M_{\ef_2}    &      & \rTo^{d} &             & M_{\ef_2}\\
\end{diagram}
If the above diagram commutes, we say that the extensions of $d$ on $M_{\ef_1}$ and $M_{\ef_2}$ {\em agree.}
In section \ref{section_agreement}, we study the analogous notions for higher derivations and conditions under which two extensions of a HD agree. Proposition \ref{extending_right} summarize that study.

The symmetric version of differential torsion theories is defined and studied in \cite{Lia_Extending}. The symmetric versions of Golan and Bland's results are shown. Also, it is shown that the symmetric Lambek, Goldie and perfect torsion theories are differential. In section \ref{section_symmetric}, we shall prove that the results from \cite{Lia_Extending} as well as the results from the previous sections of this paper hold for HDs and symmetric torsion theories (Proposition \ref{symmetric_HD_Bland_Golan}, Corollary \ref{3_symmetric_TT_are_HD} and Proposition \ref{extending_symmetric}).

Section \ref{section_questions} contains some open questions.

\section{Higher differentiability of torsion theories}
\label{section_HD}

In this section, we recall some existing results on higher differentiations and expand on some of them. First, let us consider Lemma 3.5 from \cite{Bland_paperHD}. This result is the HD version of Lemma 1.5 from \cite{Bland_paper}. Namely, Bland considers the following equivalent conditions for a Gabriel filter $\ef$ and any $n.$
\begin{itemize}
\item[($1_n$)] $\ef$ is a HD filter of order $n$.

\item[($2_n$)] For every $R$-module $M,$ every HD $\{d_n\}_{i=0}^n$ of order $n$ on $M,$ $d_i(\te(M))\subseteq\te(M)$ for all $i\leq n.$
\end{itemize}

We consider conditions:
\begin{itemize}
\item[($3'_n$)] For every $R$-module $M$, every HD $\{d_n\}_{i=0}^n$ of order $n$ on $M,$  every $x\in \te(M),$ and every $i\leq n,$ there is $K\in \ef$ such that $\delta_i(K)\subseteq\ann_r(x).$
\item[($3''_n$)] For every $R$-module $M$, every HD $\{d_n\}_{i=0}^n$ of order $n$ on $M,$  every $x\in \te(M),$ there is $K\in \ef$ such that $\delta_i(K)\subseteq\ann_r(x)$ for all $i\leq n.$
\end{itemize}
and prove that ($1_n$), ($2_n$), ($3'_n$) and ($3''_n$) are equivalent.

\begin{remark}
In \cite{Bland_paperHD}, Bland also considers the condition
\begin{itemize}
\item[($3_n$)] For every $R$-module $M$, every HD $\{d_n\}_{i=0}^n$ of order $n$ on $M,$ and every $x\in \te(M),$ there is $K\in \ef$ such that $\delta_i(K)\subseteq\ann_r(d_{n-i}(x))$ for all $i\leq n.$
\end{itemize}
as another condition claimed to be equivalent to $(1_n)$ and $(2_n).$ It is interesting to note that if the conditions $(1_n), (2_n), (3'_n)$ and $(3''_n)$ hold for $n$ that then they clearly hold for all $i\leq n.$ Condition $(3_n)$ is different in this respect.
\end{remark}

\begin{lemma} For every Gabriel filter $\ef$ and every $n$, the conditions ($1_n$), ($2_n$), ($3'_n$) and ($3''_n$) are equivalent.
\label{HD_n_Bland_lemma}
\end{lemma}
\begin{proof}
Our proof of $(1_n)\Rightarrow(2_n)$ is the same as Bland's proof of that part of Lemma 3.5 from \cite{Bland_paperHD}. Let us first note that $x\in\te(M)$ iff $\ann_r(x)$ is in $\ef.$ Thus, to prove $(2_n)$ it is sufficient to show that $\ann_r(x)\in \ef$ implies $\ann_r(d_i(x))\in \ef$ for any $i\leq n.$ We shall use induction to show that $(1_n)\Rightarrow(2_n)$. If $n=0,$ $(2_0)$ is clearly true. Assume that $(1_i)\Rightarrow(2_i)$ for all $i<n$ and assume that $(1_n)$ holds. Then $(1_i)$ holds for all $i<n$ and so we have that $(2_i)$ holds for all $i<n.$ Thus, $\ann_r(d_{i}(x))\in\ef$ for all $i<n.$ Let $I=\bigcap_{i=0}^{n-1}\ann_r(d_{i}(x)).$ Then $I$ is in $\ef.$ By $(1_n)$, for $I$ there is $J\in \ef$ such that $\delta_i(J)\subseteq I$ for all $i\leq n.$ For any $k\in J,$ $d_j(x)\delta_i(k)=0$ for all $i\leq n,$ and $j<n.$ Thus, \[d_n(x)k=d_n(xk)-\sum_{i=0}^{n-1}d_i(x)\delta_{n-i}(k)=d_n(0)-0=0\]
and so $J\subseteq\ann_r(d_n(x)).$ $J\in \ef$ implies that $\ann_r(d_n(x))\in \ef$ and so $(2_n)$ holds.

$(2_n)\Rightarrow(3'_n)$ Let $x\in \te(M).$ We shall use induction again. If $n=0,$ $(3'_0)$ is true for $K=\ann_r(x)\in \ef.$ Let us assume that $(2_i)\Rightarrow(3'_i)$ for $i<n$ and assume that $(2_n)$ holds. Then $(2_i)$ holds for all $i<n$ and so
$(3'_i)$ holds for $i<n.$ By $(2_n)$, $d_j(x)$ is in $\te(M)$ for all $j\leq n.$ Thus, by $(3'_{n-1})$, for all $d_j(x)\in\te(M),$ $j\leq n$ and all $i<n$ there are sets $K_{ij}\in \ef$ such that \[\delta_i(K_{ij})\subseteq\ann_r(d_j(x)).\] Let $K=\bigcap_{i<n, j\leq n} K_{ij}.$ $K$ is in $\ef$ and $d_j(x)\delta_i(k)=0$ for all $k\in K, j\leq n,$ and $i<n.$ Then,
\[x\delta_n(k)=d_n(xk)-\sum_{i=1}^{n}d_i(x)\delta_{n-i}(k)=0.\]
Thus, $\delta_n(K)\subseteq\ann_r(x)$ and so $(3'_n)$ holds.

$(3'_n)\Leftrightarrow(3''_n)$ Clearly $(3''_n)$ implies $(3'_n).$ Converse follows since a finite intersection of right ideals from a filter is in that filter again.

$(3''_n)\Rightarrow(1_n)$ Let $I\in \ef.$ Then $R/I$ is a torsion module and $1+I$ is a torsion element of $R/I$. By assumption, there is $K\in\ef$ such that $\delta_i(K)\subseteq \ann_r(1+I)=I$ for all $i\leq n.$ Thus, $(1_n)$ holds.
\end{proof}

\begin{lemma} If $\ef$ is a Gabriel filter corresponding to a hereditary torsion theory $\tau,$ the following conditions are equivalent.
\begin{itemize}
\item[(1)] $\ef$ is a HD filter.

\item[(2)] For every $R$-module $M,$ every $\{\delta_n\}$-HD $\{d_n\}$ on $M,$ and every $n$ \[d_n(\te(M))\subseteq\te(M).\]

\item[(3')] For every $R$-module $M$, every HD $\{\delta_n\},$ every $x\in \te(M),$ and every $n,$ there is $K\in \ef$ such that $\delta_n(K)\subseteq\ann_r(x).$

\item[(3'')] For every $R$-module $M$, every HD $\{\delta_n\},$ every $x\in \te(M),$ and every $n,$ there is $K\in \ef$ such that $\delta_i(K)\subseteq\ann_r(x)$ for all $i\leq n.$
\end{itemize}
\label{HD_Bland_lemma}
\end{lemma}
\begin{proof}
This lemma follows directly from Lemma \ref{HD_n_Bland_lemma}.
\end{proof}

In \cite{Rim1}, Rim proved that a version of part 1 Proposition \ref{Golan_Proposition} holds for rings with anti-derivations (defined on rings with involutions). In \cite{Rim2}, Rim generalized his result from \cite{Rim1} for higher anti-derivations on rings with involutions. We formulate and prove these results for higher derivations. The induction and the idea which Golan uses in his proof of Proposition \ref{Golan_Proposition} (see pages 277--279 in  \cite{Golan_paper}) is sufficient for the proof of this result. The details are below.

\begin{prop} Let $\Delta=\{\delta_n\}$ be a HD on $R$, $M$ a right $R$-module, $D=\{d_n\}$ a $\Delta$-HD on $M$ and $\tau$ a hereditary torsion theory.
\begin{enumerate}
\item If $M$ is torsion-free, then $D$ extends to a $\Delta$-HD on $M_{\ef}.$

\item If $d_n(\te(M))\subseteq \te(M)$ for all $n,$ then $D$ extends to a $\Delta$-HD on $M_{\ef}.$
\end{enumerate}
\label{HD_Golan_Prop}
\end{prop}
\begin{proof}

(1) We use induction to show that we can define a higher derivation $\{d_n\}$ on $M_{\ef}$ such that $d_n q_M=q_M d_n$ for all $n.$ For $n=0,$ the claim trivially holds. Let us assume that the claim holds for $i<n.$ Since $M$ is torsion-free, for every element $q\in M_{\ef},$ there is $I\in \ef$ such that $qI\subseteq M.$ By the induction hypothesis, $d_j(q)$ is defined for $j<n.$ Thus we can define a map $\delta_{n, q}$ on $I$ by $\delta_{n, q}(i)=d_n(qi)-\sum_{j=0}^{n-1}d_j(q)\delta_{n-j}(i).$ Note that $qi$ is in $M$ so $d_n(qi)$ is defined. The map $\delta_{n, q}$ is a homomorphism of right $R$-modules and hence it represents an element of $\dirlim\homo_R(I, M_{\ef})=(M_{\ef})_{\ef}$. But since the homomorphism
\[M_{\ef}\cong \homo_R(R, M_{\ef})\rightarrow\homo_R(I, M_{\ef})\]
is an isomorphism by Proposition 1.8, p. 198 in \cite{Stenstrom}, we can extend $\delta_{n,q}$ to $R.$ Such extension is unique by Lemma 1.9, p. 198 in \cite{Stenstrom}. Therefore, there is a unique element $q'\in M_{\ef}$ such that $q'=\delta_{n, q}(1).$  It is not hard to check that the element $q'$ is independent of the choice of ideal $I$ (using the same argument as in Golan's proof, see proof on page 278, \cite{Golan_paper}). Thus, we can define a map $d_n$ on $M_{\ef}$ by $q\mapsto q'.$ It is straightforward to check that this map is additive (using the same argument as in Golan's proof, see page 278 in \cite{Golan_paper}). To prove that $\{d_n\}$ is a higher derivation, let $r$ be any element of $R$. For $r,$ there is $J\in \ef$ such that $qrJ\subseteq M.$ $J\in\ef$ implies that $(r:J)$ is also in $\ef$ (by the definition of Gabriel filter, see \cite{Stenstrom}, page 146). Let $K=J\cap(r:J)\in \ef$ and $k\in K$

\[\delta_{n, qr}(k)-\delta_{n, q}(rk)=d_n(qrk)-\sum_{j=0}^{n-1}d_j(qr)\delta_{n-j}(k)-d_n(qrk)+\sum_{j=0}^{n-1}d_j(q)\delta_{n-j}(rk)=\]
\[-\sum_{j=0}^{n-1}\sum_{i=0}^jd_i(q)\delta_{j-i}(r)\delta_{n-j}(k)+\sum_{j=0}^{n-1}\sum_{i=0}^{n-j}d_j(q)\delta_i(r)\delta_{n-j-i}(k)\]
In the above equation, all the terms cancel except the ones with $i=n-j$ in the second sum. Thus, the above formula is equal to
\[\sum_{j=0}^{n-1}d_j(q)\delta_{n-j}(r)\delta_0(k)=\left(\sum_{j=0}^{n-1}d_j(q)\delta_{n-j}(r)\right)k\]
Hence $d_n(qr)-d_n(q)r=\sum_{j=0}^{n-1}d_j(q)\delta_{n-j}(r)$ and so $d_n(qr)=\sum_{j=0}^{n}d_j(q)\delta_{n-j}(r)$ and $\{d_n\}$ is a higher $\Delta$-derivation.

Note that $d_n(q)=\delta_{n,q}(1)$ implies that $L_{d_n(q)},$ the left multiplication by $d_n(q),$ is $\delta_{n,q}.$
On the other hand, $d_n(L_q(1))=d_n(q)=\delta_{n,q}(1)$
Thus $d_n(L_q)=\delta_{n,q}=L_{d_n(q)}.$ So, if $q$ is already in $M$, then $qR\subseteq M$ and $q_M(q)=L_q.$ In this case, $d_n(q_M(q))=d_n(L_q)=L_{d_n(q)}=q_M(d_n(q)).$

(2) Let $\overline{M}$ denote $M/\te(M)$ and $p$ denote the projection $M\rightarrow \overline{M}.$ As $d_n(\te(M))$ is contained in $\te(M),$ the additive map $\overline{d_n}:\overline{M}\rightarrow \overline{M}$ given by $\overline{d_n}p=pd_n$ is well defined for any $n$. It is easy to check that this defines a higher $\Delta$-derivation. By part (1), $\{\overline{d_n}\}$ can be extended to $\{d_n\},$ a higher $\Delta$-derivation on $M_{\ef},$ such that $d_nq_{\overline{M}}=q_{\overline{M}}\overline{d_n}$ for all $n.$ Note that $q_M=q_{\overline{M}}p$ and so $d_nq_M=d_nq_{\overline{M}}p=q_{\overline{M}}\overline{d_n} p=q_{\overline{M}}p d_n=q_Md_n.$
\end{proof}

The HD version of Bland's Theorem (Theorem \ref{Bland_Theorem}) has been proven by Bland in \cite{Bland_paperHD} (see Lemma 4.1 and Proposition 4.2 in \cite{Bland_paperHD}). For completeness, we list and prove this results.

\begin{theorem} Let $\ef$ be a Gabriel filter.
\begin{enumerate}
\item If a $\Delta$-HD on a module $M$ extends to a $\Delta$-HD on the module of quotients $M_{\ef},$ then such extension is unique.

\item The filter $\ef$ is HD if and only if for every ring HD $\Delta,$
every $\Delta$-HD on any module $M$ extends uniquely to a $\Delta$-HD on the module of quotients $M_{\ef}.$
\end{enumerate}
\label{HD_Bland_Theorem}
\end{theorem}
\begin{proof}
(1) Assume that $d'_n$ and $d''_n$ both extend $d_n$ to $M_{\ef}.$ Then $d'_n-d''_n$ is a $R$-homomorphism defined on $M_{\ef}$ that is zero on $q_M(M).$ So, it factors to an $R$-homomorphism $M_{\ef}/q_M(M)\rightarrow M_{\ef}.$ But as the module $M_{\ef}/q_M(M)$ is torsion and the module $M_{\ef}$ is torsion-free, this map has to be zero. Hence, $d'_n-d''_n=0$ on entire $M_{\ef}.$

(2) If $\ef$ is HD, the assumption of condition (2) in Proposition \ref{HD_Golan_Prop} is satisfied. Thus any HD on $M$ extends to $M_{\ef}.$ Such extension is unique by part (1) of this theorem.

Conversely, suppose that a $\Delta$-HD $\{d_n\}$ extends from $M$ to $M_{\ef}.$ For every $x\in\te(M),$
$q_M(d_n(x))=d_n(q_M(x))=d_n(0)=0.$ Thus, $d_n(x)\in \ker q_M=\te(M).$ In this case $\ef$ is HD by condition (2) of Lemma \ref{HD_Bland_lemma}.
\end{proof}

\section{Higher differentiability of Lambek, Goldie and perfect torsion theories}\label{section_Lambek_Goldie_perfect}

In this section, we prove that some frequently considered torsion theories are HD.

\begin{theorem}
\begin{enumerate}
\item The Lambek torsion theory is HD.

\item The Goldie torsion theory is HD.

\item Any perfect torsion theory is HD.
\end{enumerate}
\label{3TT_are_HD}
\end{theorem}

\begin{proof}
{\bf Part (1)}  Let $M$ be any $R$-module, $x$ any element of $M$ and $\{d_n\}$ any $\Delta$-HD. We shall prove that $\ann_r(d_n(x))$ is dense if $\ann_r(d_i(x))$ is dense for all $i<n.$ Recall that a right $R$-ideal $I$ is dense if and only if for every $r,s\in R$ such that $s\neq 0,$ there is $t\in R$ such that $st\neq 0$ and $rt\in I$. Recall also that the density of $\ann_r(x)$ implies the density of $\ann_r(xr)$ for any $r\in R$ simply because the set of all Lambek torsion elements of any module is a submodule and, as such, closed for right multiplication by elements from $R.$

Let $r,s\in R,$ and $s\neq 0.$ Since the right ideals $\ann_r(d_i(x)), i<n,$ are dense by assumption, the right ideal $I=\bigcap_{i=0}^{n-1}\ann_r(d_i(x))$ is also dense. Thus, there is $t'\in R$ such that $st'\neq 0$ and $d_i(x)rt'=0$ for all $i<n.$ Now consider the right ideal $J=\bigcap_{i<n,j\leq n}\ann_r(d_i(x)\delta_j(rt')).$ This ideal is dense as well. Thus, for 1 and $0\neq st'\in R,$ there is $t''$ such that $st't''\neq 0$ and $t''\in J$ and so $d_i(x)\delta_j(rt')t''=0$ for all $i<n,$ $j\leq n$. Let $t=t't''.$ Then $st=st't''\neq 0$ and  \[d_n(x)rt=d_n(x)rt't''=d_n(xrt')t''-\sum_{i=0}^{n-1}d_i(x)\delta_{n-i}(rt')t''=0.\]
Thus, for $r$ and $s\neq 0,$ we have found $t$ such that $rt\in \ann_r(d_n(x))$ and $st\neq0.$ Hence, $\ann_r(d_n(x))$ is dense.

{\bf Part (2)} Recall that a right ideal $I$ is essential in $R$ iff for every $r\in R,$ there is $t\in R$ such that $0\neq rt\in I.$ Let $\ef_G$ denotes the Gabriel filter of the Goldie torsion theory. Then all the essential ideals are in $\ef_G$ and
\[\ef_G=\{I\; | \mbox{ there exists } J, I\subseteq J, J\subseteq_e R\mbox{ and }(j:I)\subseteq_e R\mbox{ for all }j\in J \}.\]
For proof see Proposition 6.3, p. 148 in \cite{Stenstrom}. From this observation it is easy to see that
\[\ef_G=\{I\; | \{r\in R\; |\; (r:I) \subseteq_e R\} \subseteq_e R \}.\]

Thus, if $M$ is a right $R$-module and $x$ is an element of $M,$ $x$ is in the torsion submodule for Goldie torsion theory (i.e. $\ann_r(x)\in \ef_G$) if and only if
\[\{r\in R\; |\; \ann_r(xr) \subseteq_e R\} \subseteq_e R.\] We shall follow the notation from \cite{Lam} (see p. 255 of \cite{Lam}) and denote \[
\begin{array}{rcl}
\ann_r(x)^* & = & \{r\in R\; |\; \ann_r(xr) \subseteq_e R\}\;\;\mbox{ and }\\
\ann_r(x)^{**} & = & \{ r\in R |\, \{s\,|\, \ann_r(xrs) \subseteq_e R\}\subseteq_e R \} = \{r\in R\, |\, \ann_r(xr)^*\subseteq_e R\}.
\end{array}
\]
Using this notation, $\ann_r(x)\in\ef_G$ if and only if $\ann_r(x)^*\subseteq_e R.$ These two conditions are also equivalent to
$\ann_r(x)^{**}\subseteq_e R$ by Lemma 10 in \cite{Lia_Diff}. Note also that $J^{**}=J^{***}$ for all right ideals $J$ (see Theorem 7.28 in \cite{Lam}).

Thus, to show that condition (2) of Lemma \ref{HD_Bland_lemma} is satisfied and that the Goldie torsion theory is HD, we need to show that
\[\ann_r(d_i(x))^*\subseteq_e R\mbox{ for all }i<n\mbox{ implies that }\ann_r(d_n(x))^*\subseteq_e R.\]
So, let us assume that $\ann_r(d_i(x))^*\subseteq_e R$ for $i<n.$ Let $I=\bigcap_{i=0}^{n-1}\ann_r(d_i(x))^*.$ As $I$ is essential, for any $r\in R$ there is $t\in R$ such that $0\neq rt\in I.$ Thus $\ann_r(d_i(x)rt)\subseteq_e R$ for all $i<n.$

Let $J=\bigcap_{i=0}^{n-1}\ann_r(d_i(x)rt).$ Then $J\subseteq_e R$ so for any $s\in R$ there is $u\in R$ such that $0\neq su\in J.$ Thus $d_i(x)rtsu=0$ for all $i<n.$
\[d_n(x)rtsu=d_n(xrtsu)-\sum_{i=0}^{n-1}d_i(x)\delta_{n-i}(rtsu)\]
$d_n(xrtsu)=d_n(0)=0.$ Since $\ann_r(d_i(x))^*\subseteq_e R,$ for $i<n,$ we obtain that $\ann_r(d_i(x)\delta_{n-i}(rtsu))^*\subseteq_e R,$ $i<n$ (this is because the Goldie torsion submodule of $M$ is closed for the right multiplication with elements from $R,$ see also (2) in Lemma 11 in \cite{Lia_Diff}). Thus, the equation above gives us that $\ann_r(d_n(x)rtsu)^*\subseteq_e R.$ This proves that $\{s\;|\; \ann_r(d_n(x)rts)^* \subseteq_e R\}$ is essential. Thus, we showed that $\ann_r(d_n(x))^{***}=\{r\;|\; \{s\;|\; \ann_r(d_n(x)rs)^*\subseteq_e R\}\subseteq_e R\}$ is essential subset of $R.$ Since $\ann_r(d_n(x))^{***}=\ann_r(d_n(x))^{**}$ and $\ann_r(d_n(x))^{**}\subseteq_e R$ is equivalent to $\ann_r(d_n(x))^{*}\subseteq_e R,$ we obtain that $\ann_r(d_n(x))^{*}\subseteq_e R.$ This concludes the proof of part (2).

{\bf Part (3)} First we need an easy lemma.
\begin{lemma}
If $M$ is a right $R$-module, $N$ is an $R$-bimodule, $\{d^M_n\}$ and $\{d^N_n\}$ are two $\Delta$-HDs on $M$ and $N$ respectively, then $\{d^{M\otimes_R N}_n\},$ defined by the condition \[d^{M\otimes_R N}_n(m\otimes n)=\sum_{i=0}^n d^M_i(m)\otimes d^N_{n-i}(n),\] is a $\Delta$-HD on $M\otimes_R N.$
\label{derivation_of_tensor}
\end{lemma}
\begin{proof} Clearly all the maps $d_n^{M\otimes_R N}$ are additive.
Let $x\in M,$ $y\in N$ and $r\in R$ be arbitrary.
\begin{eqnarray}
d^{M\otimes_R N}_n((x\otimes y)r) & = & d^{M\otimes_R N}_n(x\otimes yr)\nonumber\\
& = & \sum_{i=0}^n d^M_i(x)\otimes d^N_{n-i}(nr)\nonumber\\
& = & \sum_{i=0}^n d^M_i(x)\otimes \left(\sum_{j=0}^{n-i}d^N_j(n)\delta_{n-i-j}(r)\right)\nonumber\\
& = & \sum_{i+j+k=n} d^M_i(x)\otimes d^N_j(n)\delta_{k}(r)\nonumber\\
& = & \sum_{i=0}^n \left(\sum_{j=0}^{i}d^M_i(x)\otimes d^N_{i-j}(n)\right)\delta_{n-i}(r)\nonumber\\
& = & \sum_{i=0}^n d^{M\otimes_R N}_i(x\otimes y) \delta_{n-i}(r)\nonumber
\end{eqnarray}
\end{proof}
By Theorem \ref{HD_Bland_Theorem}, it is sufficient to show that the HD $\{d_n^M\}$ extends uniquely to a HD $\{d_n^{M_{\ef}}\}$ such that $d_n^{M_{\ef}}q_M=q_M d_n^M,$ for every $n$ and every module $M.$ This is automatically fulfilled if $M$ is torsion-free by Proposition \ref{HD_Golan_Prop}.

First, let us demonstrate that we can obtain an extension $\{\delta_n^{R_{\ef}}\}$ of $\{\delta_n\}$ on $R_{\ef}.$ To prove this, it is sufficient to show that $\delta_n(\te(R))\subseteq \te(R)$ for all $n.$ Let $f: R\rightarrow R_{\ef}$ be the ring epimorphism that makes $R_{\ef}$ into a flat left $R$-module. Such map exist as $R_{\ef}$ is perfect (see Theorem 2.1, page 227 in \cite{Stenstrom}). Moreover, $f$ is the natural map $q_R.$ Thus, $\te(R)=\ker f.$ So, it is sufficient to show that $f(\delta_n(r))=0$ for all $r\in R$ such that $f(r)=0.$ Let us use induction. If $n=0$ the claim is trivially satisfied. Assume that the claim holds for all $i<n.$ Note that for $f(r)=0$ there is $m$ and $r_j\in R,$ $q_j\in R_{\ef},$ $j=1,\ldots, m$ such that $rr_j=0$ and $\sum_j f(r_j)q_j=1$ (see part c in Theorem 2.1, page 227 in \cite{Stenstrom}). Thus,
\[\begin{array}{rcll}
f(\delta_n(r)) & = & f(\delta_n(r)) 1 = f(\delta_n(r))\sum_j f(r_j)q_j & (\mbox{by above})\\
 & = & \sum_j f\left(\delta_n(r)r_j\right)q_j & (f\mbox{ is a ring hom.})\\
& = & \sum_j f\left(\delta_n(rr_j)-\sum_{i=0}^{n-1}\delta_i(r)\delta_{n-i}(r_j)\right)q_j & (\{\delta_n\}\mbox{ is a HD})\\
& = & \sum_j f\left(-\sum_{i=0}^{n-1}\delta_i(r)\delta_{n-i}(r_j)\right)q_j & (rr_j=0\mbox{ for all }j)\\
& = & -\sum_j \sum_{i=0}^{n-1}f(\delta_i(r))f(\delta_{n-i}(r_j))q_j  & (f\mbox{ is a ring hom.})\\
& = & 0 & (f(\delta_i(r))=0\mbox{ for }i<n)\\
\end{array}\]

As $\ef$ is perfect, the unique map $Q_M: M\otimes_R R_{\ef}\rightarrow
M_{\ef}$ such that $q_M = Q_M i_M$ is an isomorphism for every module $M.$
Define \[d_n^{M_{\ef}}=Q_M d_n^{M\otimes_R R_{\ef}}Q_M^{-1}\]
where map $d_n^{M\otimes_R R_{\ef}}$ is the map from Lemma \ref{derivation_of_tensor} defined via $\{d_n^M\}$ and $\{\delta_n^{R_{\ef}}\}.$

Clearly, the maps $d_n^{M_{\ef}},$ $n\in \omega$ are additive and a straightforward calculation shows that $\{d_n^{M_{\ef}}\}$ is a higher derivation. We show that the following diagram commutes
\begin{diagram}
M & \rTo^{i_M} & M\otimes_R R_{\ef} & \rTo^{Q_M} & M_{\ef}\\
\dTo_{d_n^M} & & \dTo_{d_n^{M\otimes_R R_{\ef}}} & & \dTo_{d_n^{M_{\ef}}} \\
M & \rTo^{i_M} & M\otimes_R R_{\ef} & \rTo^{Q_M} & M_{\ef}
\end{diagram}
As $\delta_n(1)=0,$ $\delta_n^{R_{\ef}}(1)=0.$  Thus, if $m\in M$ is arbitrary, $d_n^{M\otimes_R R_{\ef}}i_M(m)=d_n^{M\otimes_R R_{\ef}}(m\otimes 1)=d_n^{M}(m)\otimes 1+0=i_M(d_n^M(m)).$ So, the maps in the first square commute. The maps in the second diagram commute by the definition of $d_n^{M_{\ef}}.$
This gives us

\[\begin{array}{rcll}
d_n^{M_{\ef}}q_M & = & Q_M d_n^{M\otimes_R R_{\ef}}Q_M^{-1}q_M & (\mbox{definition of }d_n^{M_{\ef}})\\
& = & Q_M d_n^{M\otimes_R R_{\ef}}i_M  & (\mbox{as }q_M = Q_M i_M)\\
& = & Q_M i_M d_n^M & (\mbox{by the above diagram})\\
& = & q_M d_n^M  & (\mbox{as }q_M = Q_M i_M).
\end{array}\]
Finally, $d_n^{M_{\ef}}$ is unique by Theorem \ref{HD_Bland_Theorem}.
This finishes the proof of the Theorem \ref{3TT_are_HD}.
\end{proof}

\section{Agreement of different extensions of a HD}
\label{section_agreement}

In this section we study the conditions under which extensions of a HD on two different modules of quotients agree. We prove that the results from \cite{Lia_Extending} generalize to higher derivations as well. First, let  $\ef_1$ and $\ef_2$ be two filters such that $\ef_1\subseteq \ef_2.$ Let $M$ be a right $R$-module and $\{d_n\}$ a $\Delta$-HDs defined on $M.$ If $\{d_n\}$ extends to $M_{\ef_1}$ and $M_{\ef_2}$ in such a way that the following diagram commutes for every $n,$ then we say that the extensions of $\{d_n\}$ on $M_{\ef_1}$ and $M_{\ef_2}$ {\em agree.}
\begin{diagram}
  &             & M_{\ef_1}    &      & \rTo^{d_n} &             & M_{\ef_1}\\
  & \ruTo^{q_1} & \vLine       &      &            & \ruTo^{q_1} &    \\
M &             & \rTo^{d_n}     &      & M          &             & \dTo_{q_{12}}\\
  & \rdTo^{q_2} & \dTo_{q_{12}}&      &            & \rdTo^{q_2} &     \\
  &             & M_{\ef_2}    &      & \rTo^{d_n} &             & M_{\ef_2}\\
\end{diagram}

\begin{prop}
Suppose that a Gabriel filter $\ef_1$ is contained in a Gabriel filter $\ef_2$ and that $M$ is a right $R$-module with a $\Delta$-HD $\{d_n\}$. If $d_n$ can be extended to $M_{\ef_1}$ for every $n$ and either
\begin{itemize}
\item[i)] $d_n$ can be extended from $M_{\ef_1}$ to $M_{\ef_2}$ for every $n,$  or

\item[ii)] $d_n$ can be extended from $M$ to $M_{\ef_2}$ for every $n,$
\end{itemize}
then the extensions of $\{d_n\}$ to $M_{\ef_1}$ and $M_{\ef_2}$ agree.

As a consequence, if $\ef_1$ is HD and either
\begin{itemize}
\item[i)] $M_{\ef_1}$ is $\tau_2$-torsion-free, or

\item[ii)] $\ef_2$ is a HD filter,
\end{itemize}
then any HD on $M$ extends both to $M_{\ef_1}$ and $M_{\ef_2}$ in such a way that the extensions on $M_{\ef_1}$ and $M_{\ef_2}$ agree.
\label{extending_right}
\end{prop}
\begin{proof}
Since $\ef_1\subseteq\ef_2,$ $q_{12}q_1=q_2.$ By assumption that $\{d_n\}$ can be extended to a higher derivation on $M_{\ef_1},$ we have that $d_n q_1=q_1 d_n.$

In case i), we have that $d_nq_{12}=q_{12} d_n$ and need to prove that $d_nq_2=q_2 d_n.$ This is the case because
\[d_nq_2=d_n q_{12} q_1 = q_{12} d_n q_1 = q_{12} q_1 d_n = q_2 d_n.\]

In case ii), we have that $d_nq_2=q_2 d_n$ and need to prove that $d_nq_{12}=q_{12} d_n.$ Clearly, this is the case for $n=0.$ Let us assume that this is the case for $i<n$ and let us show that the claim is true for $n$. Let $q\in M_{\ef_1}.$ Then there is a right ideal $I$ in $\ef_1$ such that $qI\subseteq M/\te_1(M).$ In this case $q_{12}(q)I\subseteq M/\te_2(M)$ by definition of map $q_{12}.$ Note also that if $qr=m+\te_1(M)$ for some $r\in I$ and $m\in M,$ then $q_{12}(q)r=m+\te_2(M)$ since
\[L_{q_{12}(q)r}=q_{12}(L_{qr})=q_{12}(L_{m+\te_1(M)})=q_{12}(q_1(m))=q_2(m)=L_{m+\te_2(M)}.\]

As a higher derivation on $M$ can be extended to $M_{\ef_1},$ $d_n(\te_1(M))\subseteq\te_1(M)$ for all $n$ by Theorem \ref{HD_Bland_Theorem}. Thus $\{d_n\}$ defines a higher derivation $\{\overline{d_n}\}$ on $M/\te_1(M)$ such that $\overline{d_n}(m+\te_1(M))=d_n(m)+\te_1(M)$ for all $n.$ The extension of $d_n$ on $M_{\ef_1}$ coincides with the extension of $\overline{d_n}$ from $M/\te_1(M)$ to $M_{\ef_1}$ (see the proof of part (2) of Proposition \ref{HD_Golan_Prop}). Thus $q_1(d_n(m))=q_1(\overline{d_n}(m+\te_1(M)))=d_1(q_n(m)).$
Similarly, $\{d_n\}$ defines a higher derivation $\{\overline{d_n}\}$ on $M/\te_2(M)$ such that $q_2(d_n(m))=q_2(\overline{d_n}(m+\te_2(M)))=d_n(q_2(m)).$

Thus, the extension $d_n$ on $M_{\ef_1}$ is defined such that \[d_n(q)r=d_n(qr)- \sum_{i=0}^n d_i(q)\delta_{n-i}(r) =L_{\overline{d_n}(qr)}-\sum_{i=0}^n d_i(q)\delta_{n-i}(r)\] for all $r\in I.$
As $I$ is in $\ef_2$ also and $q_{12}(q)I\subseteq M/\te_2(M)$, we have that \[d_n(q_{12}(q))r= d_n(q_{12}(q)r)-\sum_{i=0}^n d_iq_{12}(q)\delta_{n-i}(r)=L_{\overline{d_n}(q_{12}(q)r)}-\sum_{i=0}^n d_iq_{12}(q)\delta_{n-i}(r).\]

Let us assume now that $d_i (q_{12}(q))r=q_{12}(d_i(q))r$ for $i<n$ and show that $d_n (q_{12}(q))r=q_{12}(d_n(q))r.$
\[
\begin{array}{rcll}
d_n (q_{12}(q))r & = & L_{\overline{d_n}(q_{12}(q)r)}-\sum_{i=0}^n d_iq_{12}(q)\delta_{n-i}(r) & (\mbox{see above})\\
& =&  L_{\overline{d_n}(m+\te_2(M))}-\sum_{i=0}^n d_iq_{12}(q)\delta_{n-i}(r) & (\mbox{def. of }m)  \\
& =&  d_n(L_{m+\te_2(M)})-\sum_{i=0}^n d_iq_{12}(q)\delta_{n-i}(r) & (\mbox{see above}
)\\
& =&  d_nq_2(m)-\sum_{i=0}^n d_iq_{12}(q)\delta_{n-i}(r) & (\mbox{by def. of }q_2)\\
& =& q_2d_n(m)-\sum_{i=0}^n d_iq_{12}(q)\delta_{n-i}(r) & (d_nq_2 =q_2 d_n)\\
& =& q_{12}q_1d_n(m)-\sum_{i=0}^n d_iq_{12}(q)\delta_{n-i}(r) & (q_2=q_{12}q_1)\\
& =&  q_{12}d_nq_1(m)-\sum_{i=0}^n d_iq_{12}(q)\delta_{n-i}(r) & (d_nq_1 =q_1 d_n)\\
& =&  q_{12}(d_n(L_{m+\te_1(m)}))-\sum_{i=0}^n d_iq_{12}(q)\delta_{n-i}(r) & (\mbox{by def. of }q_1)\\
& =&  q_{12}(L_{\overline{d_n}(m+\te_1(m))})-\sum_{i=0}^n d_iq_{12}(q)\delta_{n-i}(r) & (\mbox{see above}
)\\
& =&  q_{12}(L_{\overline{d_n}(qr)})-\sum_{i=0}^n d_iq_{12}(q)\delta_{n-i}(r) &  (\mbox{def. of }m)\\
& =&  q_{12}(L_{\overline{d_n}(qr)})-\sum_{i=0}^n q_{12}d_i(q)\delta_{n-i}(r)) & (\mbox{by induction})\\
& = & q_{12}\left(L_{\overline{d_n}(qr)}-\sum_{i=0}^n d_i(q)\delta_{n-i}(r)\right) & (q_{12}\mbox{ is an }R\mbox{-map})\\
& =&  q_{12}(d_n(q)r) & (\mbox{see above})\\
& =&  q_{12}(d_n(q))r & (q_{12}\mbox{ is an }R\mbox{-map})\\
\end{array}
\]

Thus, the left multiplication with $d_n (q_{12}(q))-q_{12}(d_n(q))$ defines the zero $R$-map $I\rightarrow M_{\ef_2}.$ As $I\in\ef_2,$ this map extends to a $R$ map $f: R\rightarrow M_{\ef_2}$ (see Proposition 1.8 p. 198 in \cite{Stenstrom}). Since $I\subseteq \ker f,$ $f$ factors to a map $R/I\rightarrow M_{\ef_2}.$ But this map has to be zero as $R/I$ is torsion and $M_{\ef_2}$ is torsion-free in $\tau_2.$  Thus, $f$ is zero and so $d_n (q_{12}(q))=q_{12}(d_n(q))$ for every $q\in M_{\ef_1}.$

The second part of the proposition follows directly from the first part. If $M_{\ef_1}$ is $\tau_2$-torsion-free, then $\{d_n\}$ can be extended from $M_{\ef_1}$ to $M_{\ef_2}$ by part (1) of Proposition \ref{HD_Golan_Prop}. The extensions agree by the first part of this proposition. If $\ef_2$ is HD, then $\{d_n\}$ can be extended from $M$ to $M_{\ef_2}.$ In this case part ii)  of the first part of the proposition is satisfied and thus the extensions on $M_{\ef_1}$ and $M_{\ef_2}$ agree.
\end{proof}

A corollary of Theorem \ref{3TT_are_HD} and Proposition \ref{extending_right} is that an extension of any HD to a module (or ring) of quotients with respect to any HD, hereditary and faithful torsion theory agrees with the extension with respect to the Lambek torsion theory. In particular, a HD on $R$ extends to the maximal right ring of quotients $\Qmax(R)$ in such a way that its extension agrees with the extension on the total right ring of quotients $\Qtot(R)$ (the maximal perfect right ring of quotients). Also, the extension on the right ring of quotients with respect to the Goldie torsion theory agrees with the extensions on $\Qmax(R)$ and $\Qtot(R).$ If $R$ is right Ore, these extensions agree with the extension on the classical right ring of quotients $\Qcl(R). $

\section{Higher differentiability of symmetric torsion theories}
\label{section_symmetric}

In this section, we shall prove the symmetric versions of results from the previous sections.

If $R$ and $S$ are two rings, $\ef_l$ a Gabriel filter of left $R$-ideals and $\ef_r$ a Gabriel filters of right $S$-ideals, define $_{\ef_l}\ef_{\ef_r}$ as the set of right ideals of $S\otimes_{\Zset}
R^{op}$ containing ideals of the form $J\otimes R^{op}+S\otimes I$ where $I\in\ef_l$ and $J\in \ef_r.$ This defines a Gabriel filter (for details see \cite{Ortega_paper}). We shorten the notation $_{\ef_l}\ef_{\ef_r}$ to
$_{l}\ef_{r}$ when there is no confusion about the Gabriel filters used. The filter $_{l}\ef_{r}$ is called the {\em symmetric filter induced by $\ef_l$ and $\ef_r$}. If $\tau_l$ and $\tau_r$ are the torsion theories corresponding to filters $\ef_l$ and $\ef_r$ respectively, we call the torsion theory $_l\tau_r$ corresponding to $_{l}\ef_{r}$ {\em the torsion theory induced by $\tau_l$ and $\tau_r.$}

If $M$ is an $R$-$S$-bimodule, $\te_l(M),$ $\te_r(M)$ and $_l\te_r(M)$ torsion submodules of $M$ for $\tau_l,$ $\tau_r$ and  $_l\tau_r$ respectively, then \[_l\te_r(M)=\te_l(M)\cap\te_r(M).\] For details see \cite{Golan}, I, ch. 2, Proposition 2.5. Thus the torsion theory on $S\otimes_{\Zset}
R^{op}$ corresponding to filter $_l\ef_r$ of right $S\otimes_{\Zset}
R^{op}$-ideals is exactly the torsion theory of $R$-$S$-bimodules with the torsion class $\te_l\cap\te_r.$

In \cite{Ortega_paper}, Ortega defines the {\em symmetric module of quotients} $_{\ef_l}M_{\ef_r}$ of $M$ with respect to $_l\ef_r$ to be
\[_{\ef_l}M_{\ef_r}=\dirlim_{K\in _l\ef_r}\;  \homo(K, \frac{M}{_l\te_r(M)})\]
where the homomorphisms in the formula are $S\otimes R^{op}$ homomorphisms (equivalently $R$-$S$-bimodule homomorphisms).
We shorten the notation $_{\ef_l}M_{\ef_r}$ to $_lM_r.$ Just as in the right-sided case, there is a left exact functor $q_M$ mapping $M$ to the symmetric module of quotients $_lM_r$ obtained by composing the projection $M\rightarrow M/\,_l\te_r(M)$ with the injection $M/\,_l\te_r(M)\rightarrow\, _lM_r$ (see Lemma 4 in \cite{Lia_Sym}).

If $R=S,$ the module $_lR_r$ has a ring structure (Lemma 1.5 in \cite{Ortega_paper}). The ring $_lR_r$ is called {\em the symmetric ring of quotients} with respect to the torsion theory $_l\tau_r.$

Note that every derivation $\delta$ on $R$ determines a derivation on $R\otimes_{\Zset}R^{op}$ given by \[\overline{\delta}(r\otimes s)=\delta(r)\otimes s+r\otimes \delta(s).\] Similarly, every HD $\Delta$ on $R$ determines a HD $\overline{\Delta}$ on $R\otimes_{\Zset}R^{op}$ given by
\[\overline{\delta_n}(r\otimes s)=\sum_{i=0}^n \delta_i(r)\otimes \delta_{n-i}(s).\]

If $M$ is an $R$-bimodule, and $\delta$ a derivation on $R$, an additive map $d: M\rightarrow M$ is a {\em $\delta$-derivation} if
\[d(xr)=d(x)r+x\delta(r)\mbox{ and }d(rx)=\delta(r)x+rd(x)\]
for all $x\in M$ and $r\in R.$ It is straightforward to check that $d$ is a $\overline{\delta}$-derivation on $M$ considered as a right $R\otimes_{\Zset}R^{op}$-module. Conversely, every $\overline{\delta}$-derivation of a right $R\otimes_{\Zset}R^{op}$-module $M$ determines a $\delta$-derivation on $M$ considered as an $R$-bimodule.

This generalizes to HDs as well. If $M$ is an $R$-bimodule, and $\Delta$ a HD on $R$, we shall say that an indexed family of additive maps $\{d_n\}$ defined on $M$ is a {\em $\Delta$-HD} if $d_0$ is an identity,
\[d_n(xr)=\sum_{i=0}^n\delta_i(x)\delta_{n-i}(r)\mbox{ and }d_n(rx)=\sum_{i=0}^n\delta_i(r)\delta_{n-i}(x)\]
for all $x\in M$ and $r\in R.$ It is straightforward to check that $\{d_n\}$ is a $\overline{\Delta}$-HD on $M$ considered as a right $R\otimes_{\Zset}R^{op}$-module. Conversely, every $\overline{\Delta}$-HD on a right $R\otimes_{\Zset}R^{op}$-module $M$ determines a $\Delta$-HD on $M$ considered as an $R$-bimodule. Specifically, every HD $\Delta$ on $R$ is a $\overline{\Delta}$-HD on $R$ considered as a right $R\otimes_{\Zset}R^{op}$-module. Conversely, every HD $\overline{\Delta}$ on $R\otimes_{\Zset}R^{op}$ is a $\Delta$-HD on $R\otimes_{\Zset}R^{op}$ considered as an $R$-bimodule.

In \cite{Lia_Extending}, a symmetric filter $_l\ef_r$ induced by a left Gabriel filter $\ef_l$ and a right Gabriel filter $\ef_r$ is said to be {\em differential} if for every $I\in\, _l\ef_r$ there is $J\in\, _l\ef_r$ such that $\overline{\delta}(J)\subseteq I$ for all $R\otimes_{\Zset}R^{op}$ derivations $\overline{\delta}.$ If we consider the right $R\otimes_{\Zset}R^{op}$-ideals $I$ and $J$ as $R$-bimodules, the condition $\overline{\delta}(J)\subseteq I$ is equivalent with $\delta(J)\subseteq I$ (for more details see section 3 in \cite{Lia_Extending}). The hereditary torsion theory determined by $_l\ef_r$ is said to be {\em differential} in this case.

\begin{definition}
A symmetric filter $_l\ef_r$ is a {\em HD symmetric filter} if for every $R\otimes_{\Zset}R^{op}$ HD $\{\overline{\delta_n}\},$ every $I\in\, _l\ef_r,$ and every $n,$ there is $J\in\, _l\ef_r$ such that $\overline{\delta_n}(J)\subseteq I.$ The hereditary torsion theory determined by $_l\ef_r$ is said to be {\em higher differential} (HD for short) in this case.
\end{definition}

If we consider the right $R\otimes_{\Zset}R^{op}$-ideals $I$ and $J$ as $R$-bimodules, the condition $\overline{\delta_n}(J)\subseteq I$ is equivalent with $\delta_n(J)\subseteq I$ (by observations on symmetric HDs above).

\begin{definition}
Let $\Delta$ be a HD on $R$, $M$ an $R$-bimodule, $D=\{d_n\}$ a $\Delta$-HD on $M$ and $_l\tau_r$ a symmetric torsion theory. If $D$ is such that every $d_n$ extends to the module of quotients $_lM_r$ such that $d_n q_M=q_M d_n$ for all $n,$ then we shall say that $D$ {\em extends} to a $\Delta$-HD on $_lM_r.$
\end{definition}

The symmetric versions of Bland and Golan's results (Proposition \ref{Golan_Proposition} and Theorem \ref{Bland_Theorem}) on derivations have been proven in \cite{Lia_Extending} (Proposition 3 and Theorem 2 in \cite{Lia_Extending}). From the proof of those two results and from Lemma \ref{HD_Bland_lemma}, Proposition \ref{HD_Golan_Prop} and Theorem \ref{HD_Bland_Theorem} from section \ref{section_HD} of this paper, we arrive to the symmetric versions of Bland and Golan's results generalized to HDs. In what follows, we use $_l\ann_r(x)$ to denote $\{t\in R\otimes R^{op}| xt=0\}$ for an element $x$ of an $R$-bimodule.

\begin{prop} Let $_l\ef_r$ be a symmetric Gabriel filter induced by $\ef_l$ and $\ef_r,$ corresponding to the symmetric torsion theory $_l\tau_r.$
\begin{itemize}
\item[(i)] {\em Lemma \ref{HD_Bland_lemma} for symmetric filters.} The following conditions are equivalent.
\begin{itemize}
\item[(1)] $_l\ef_r$ is a HD symmetric filter.

\item[(2)] For every $R$-bimodule $M,$ every $\Delta$-HD $\{d_n\}$ on $M,$ and every $n,$ \[d_n(\,_l\te_r(M))\subseteq\, _l\te_r(M).\]

\item[(3')] For every $R$-bimodule $M,$ every HD $\Delta,$ every $x\in\,_l\te_r(M),$ and every $n,$ there is $K\in\,_l\ef_r$ such that $\overline{\delta_n}(K)\subseteq\, _l\ann_r(x).$

\item[(3'')] For every $R$-bimodule $M,$ every HD $\Delta,$ every $x\in\,_l\te_r(M),$ and every $n,$ there is $K\in\,_l\ef_r$ such that $\overline{\delta_i}(K)\subseteq\, _l\ann_r(x)$ for all $i\leq n.$
\end{itemize}
Moreover, if $\ef_l$ and $\ef_r$ are HD filters, then $_l\ef_r$ is a HD symmetric filter.

\item[(ii)] {\em Golan's Proposition for HDs and symmetric filters.}
If $\Delta$ is a ring HD, $M$ an $R$-bimodule, $D=\{d_n\}$ a $\Delta$-HD on $M,$ then
\begin{enumerate}
\item If $M$ is torsion-free, then $D$ extends to a $\Delta$-HD on the symmetric module of quotients $_lM_r.$

\item If $d_n(\, _l\te_r(M))\subseteq\, _l\te_r(M)$ for all $n,$ then $D$ extends to a $\Delta$-HD on the symmetric module of quotients $_lM_r.$
\end{enumerate}

\item[(iii)] {\em Bland's Theorem for HDs and symmetric filters.}
\begin{enumerate}
\item If a $\Delta$-HD on a bimodule $M$ extends to a $\Delta$-HD on the symmetric module of quotients $_lM_r,$ then such extension is unique.

\item The filter $_l\ef_r$ is a HD symmetric filter if and only if for every ring HD $\Delta,$
every $\Delta$-HD on any bimodule $M$ extends uniquely to a $\Delta$-HD on the symmetric module of quotients $_lM_r.$
\end{enumerate}
\end{itemize}
\label{symmetric_HD_Bland_Golan}
\end{prop}

\begin{proof}
The proof of the equivalence of four conditions in part (i) follows the proof of Proposition 3 in \cite{Lia_Extending} using induction just as in the proof of Lemma \ref{HD_n_Bland_lemma}.

To prove the last sentence of part (i), let us assume that $\ef_l$ and $\ef_r$ are HD. Let $M$ be an $R$-bimodule, $\Delta$ a HD on $R$ and $D$ a $\Delta$-derivation on $M.$ As $\ef_l$ and $\ef_r$ are HD,
$d_n(\te_l(M))\subseteq \te_l(M)$ and $d_n(\te_r(M))\subseteq \te_r(M)$ for any $n$ by Lemma \ref{HD_Bland_lemma}. Thus
$d_n(_l\te_r(M))=d_n(\te_l(M)\cap \te_r(M)))\subseteq d_n(\te_l(M))\cap d_n(\te_r(M))$ is contained in $\te_l(M)\cap\te_r(M)=\,_l\te_r(M).$ Thus condition (2) of Lemma \ref{HD_n_Bland_lemma} holds and so $_l\ef_r$ is HD.

The proof of (1) of part (ii) follows the proof of (1) in Theorem 2 in \cite{Lia_Extending} as follows. Consider $M$ as a right $R\otimes_{\Zset}R^{op}$-module. Note that $_lM_r$ is a right $R\otimes_{\Zset}R^{op}$-module of quotients with respect to $_l\tau_r$ considered as the torsion theory of right $R\otimes_{\Zset}R^{op}$-modules. Thus,
$D$ extends to a $\overline{\Delta}$-HD on $_lM_r$ by (1) of Proposition \ref{HD_Golan_Prop}. But then such extension is a $\Delta$-HD on $_lM_r.$

The proof of (2) of part (ii) follows the proof of (2) in Proposition \ref{HD_Golan_Prop}.

The proof of (1) and (2) of part (iii) follow the proofs of (1) and (2) of Theorem \ref{HD_Bland_Theorem}. In proof of (2) the part (i) is used instead of Lemma \ref{HD_Bland_lemma}.
\end{proof}

The symmetric version of the Lambek torsion theory is the symmetric torsion theory induced by the Gabriel filters of all dense right and left $R$-ideals respectively. It is called the symmetric Lambek theory and it is the largest symmetric hereditary and faithful torsion theory. Its symmetric ring of quotients is called the maximal symmetric ring of quotients and is denoted it by $\Qsimmax(R)$ in \cite{Lia_Sym} and \cite{Lia_Extending}.

The symmetric Goldie torsion theory is the theory induced by the torsion theories of left and right $R$-modules whose torsion-free classes consist of left and right nonsingular modules respectively.

In \cite{Lia_Sym}, the symmetric versions of right perfect rings of quotients and the total right ring of quotients are defined and studied. A ring extension $S$ of $R$ with an embedding $f: R\rightarrow S$ is a perfect symmetric ring of quotients if the family of left ideals $\ef_l=\{I | Sf(I)=S\}$ is a
left Gabriel filter, the family of right ideals $\ef_r=\{J | f(J)S=S\}$ is a right
Gabriel filter and there is a ring isomorphism $g: S\cong\; _{\ef_l}R_{\ef_r}$ such that
$g\circ f$ is the canonical map $q_R: R\rightarrow\; _{\ef_l}R_{\ef_r}.$ For the equivalent versions of this definition see Theorem 11 of \cite{Lia_Sym}. In this setup, $\ef_r$ is a perfect right, $\ef_l$ a perfect left and $_l\ef_r,$ the filter induced by $\ef_l$ and $\ef_r,$ is called a perfect symmetric filter. The torsion theory $_l\ef_r$ determines is called a perfect symmetric torsion theory.

Symmetrizing the condition for the perfect right torsion theories $M_{\ef}\cong M\otimes_R R_{\ef}$, we arrive to the condition describing the perfect symmetric torsion theories and symmetric filters:
\[_lR_r \otimes_R M\otimes_R\, _lR_r\cong \;_lM_r\]
for all $R$-bimodules $M.$ To check out the proof of this fact as well as the list of equivalent conditions, see Theorem 12 in \cite{Lia_Sym}.

As in the one-sided case, the following result holds.
\begin{corollary}
\begin{enumerate}
\item The symmetric Lambek torsion theory is HD.

\item The symmetric Goldie torsion theory is HD.

\item Any perfect symmetric torsion theory is HD.
\end{enumerate}
\label{3_symmetric_TT_are_HD}
\end{corollary}
\begin{proof}
Since the Lambek and Goldie torsion theories of right (and left) modules are HD by Theorem \ref{3TT_are_HD} and since the symmetric Lambek and Goldie torsion theories are induced by one-sided Lambek and Goldie theories respectively, (1) and (2) directly follow from the last sentence in part (i) of Proposition \ref{symmetric_HD_Bland_Golan}.

As pointed out in paragraphs preceding the corollary, a perfect symmetric filter $_l\ef_r$ is induced by filters $\ef_l$ and $\ef_r$ that are perfect filters. So, $\ef_l$ and $\ef_r$ are HD by Theorem \ref{3TT_are_HD}. Then $_l\ef_r$ is also HD by part (i) of Proposition \ref{symmetric_HD_Bland_Golan}.
\end{proof}

The results on extending a HD (summarized in Proposition \ref{extending_right}) also hold for symmetric filters.
Let $\ef^1_l$ and $\ef^2_l$ be left and $\ef^1_r$ and $\ef^2_r$ be right Gabriel filters that induce the symmetric filters $_l\ef_r^1$ and $_l\ef_r^2.$ Let $_l\tau_r^i$ denote the corresponding torsion theories. Suppose that $_l\ef_r^1$ is contained in $_l\ef_r^2.$ If $M$ is an $R$-bimodule, let $q_i$ denote the left exact functors mapping $M$ to the rings of quotients $_lM_r^i$ with respect to $_l\ef_r^i$ for $i=1$ and $2.$ Note that we have the mapping $q_{12}:\,_lM_r^1\rightarrow\,_lM_r^2$ induced by the inclusion $_l\ef_r^1\subseteq\,_l\ef_r^2$ such that $q_{12}q_1=q_2$ just as in the right-sided case.

If $D=\{d_n\}$ is a $\Delta$-HD defined on $M$ and if $d_n$ extends to $_lM_r^1$ and $_lM_r^2$ in such a way that the following diagram commutes for every $n,$ we say that the extensions of $D$ on $_lM_r^1$ and $_lM_r^1$ {\em agree.}
\begin{diagram}
  &             & _lM_r^1      &     & \rTo^{d_n} &             & _lM_r^1\\
  & \ruTo^{q_1} & \vLine       &     &            & \ruTo^{q_1} &    \\
M &             & \rTo^{d_n}     &     & M          &             & \dTo_{q_{12}}\\
  & \rdTo^{q_2} & \dTo_{q_{12}}&     &            & \rdTo^{q_2} &     \\
  &             & _lM_r^2      &     & \rTo^{d_n} &             & _lM_r^2\\
\end{diagram}

\begin{prop} Suppose that a symmetric filter $_l\ef_r^1$ is contained in a symmetric filter $_l\ef_r^2$ and that $\{d_n\}$ is a $\Delta$-HD on a bimodule $M.$ If $d_n$ can be extended to $_lM_r^1$ for every $n$ and either
\begin{itemize}
\item[i)] $d_n$ can be extended from $_lM_r^1$ to $_lM_r^2$ for every $n$  or

\item[ii)] $d_n$ can be extended from $M$ to $_lM_r^2$ for every $n,$
\end{itemize}
then the extensions of $\{d_n\}$ to $_lM_r^1$ and $_lM_r^2$ agree.

As a consequence, if $_l\ef_r^1$ is HD and either
\begin{itemize}
\item[i)] $_lM_r^1$ is $_l\tau_r^2$-torsion-free, or

\item[ii)] $_l\ef_r^2$ is HD,
\end{itemize}
then any HD on $M$ extends both to $_lM_r^1$ and $_lM_r^2$ in such a way that the extensions agree.
\label{extending_symmetric}
\end{prop}
\begin{proof}
The proof of the first part follows from the proof of Proposition \ref{extending_right} exactly using part (iii) of Proposition \ref{symmetric_HD_Bland_Golan} instead of Theorem \ref{HD_Bland_Theorem}.

The second part follows directly from the first part similarly as the second part of Proposition \ref{extending_right} follows from the first part (just use part (ii) of Proposition \ref{symmetric_HD_Bland_Golan} instead of Proposition \ref{HD_Golan_Prop}).
\end{proof}

By Theorem 14 in \cite{Lia_Sym}, every ring has the largest perfect symmetric ring of quotients -- the total symmetric ring of quotients. We denote it by $\Qsimtot(R)$ as in \cite{Lia_Sym}. Analogously to the right-sided case, $R\subseteq \Qsimtot(R)\subseteq \Qsimmax(R).$ If $R$ is Ore, $R\subseteq\Qlrcl(R)\subseteq \Qsimtot(R)\subseteq \Qsimmax(R)$ (see Corollary 16 in \cite{Lia_Sym}).

From Corollary \ref{3_symmetric_TT_are_HD} and Proposition \ref{extending_symmetric} it follows that an extension of any HD to a module (or ring) of quotients with respect to any HD, hereditary and faithful symmetric torsion theory agrees with the extension with respect to the symmetric Lambek torsion theory. In particular, a HD on $R$ extends to the maximal symmetric ring of quotients $\Qsimmax(R)$ in such a way that its extension agrees with the extension on the total symmetric ring of quotients $\Qsimtot(R).$ Also, the extension on the symmetric ring of quotients with respect to the symmetric Goldie torsion theory agrees with the extensions on $\Qsimmax(R)$ and $\Qsimtot(R).$ If $R$ is Ore, these extensions agree with the extension on the classical ring of quotients $\Qlrcl(R).$

\section{Questions and open problems}
\label{section_questions}

We list some of the open questions. The first two have also been raised in \cite{Lia_Extending}.
\begin{enumerate}
\item Is every hereditary torsion theory differential? All the hereditary torsion theories that are most frequently in use (classical, perfect, Goldie, Lambek) are differential. Is there a non-example?

\item Let us suppose that $\ef_1$ and $\ef_2$ are right (or symmetric) filters with $\ef_1\subseteq\ef_2.$ Let $M$ be a right (or bi) $R$-module with a derivation $d.$ If $\ef_2$ is differential, when can the extension of $d$ to the module of quotient with respect to $\ef_2$ be restricted to a derivation on the module of quotients with respect to $\ef_1$ in such a way that the derivations agree (in the sense used throughout this paper and \cite{Lia_Extending})? This question might be related to the first one since every hereditary and faithful torsion is contained in the Lambek torsion theory that is differential.

\item Find an example of a torsion theory that is differential but not higher differential.
\end{enumerate}

\section*{Acknowledgments}

After studying the extendibility of ring derivations, the question of the extendibility of {\em higher} derivations was first brought to the author's attention by Wagner de Oliveira Cortes. The author is deeply grateful to Wagner de Oliveira Cortes for that.

\end{document}